\documentclass[12pt]{article}

\usepackage{amsfonts,amsmath,graphicx}
\usepackage[final]{epsfig}
\usepackage{graphics}
\usepackage{amsmath}
\usepackage{amsfonts}
\usepackage{latexsym}
\usepackage{amssymb}
\usepackage{graphicx}
\usepackage{epstopdf}
\newtheorem{lemma}{Lemma}[section]

\newtheorem{remark}[lemma]{Remark}

\newtheorem{theorem}{Theorem}

\newtheorem{corollary}[theorem]{Corollary}

\newtheorem{conjecture}[theorem]{Conjecture}

\newcommand{\proofend}{$\Box$\bigskip}

\newcommand{\Q}{{\mathbb Q}}
\newcommand{\R}{{\mathbb R}}
\newcommand{\Z}{{\mathbb Z}}
\newcommand{\RP}{{\mathbb{RP}}}

\def\proof{\paragraph{Proof.}}

\begin{document}
\date{\today}

\title {Periodic  trajectories in the regular pentagon}
\author{Diana Davis,\thanks{
Department of Mathematics,
Brown University, Providence, RI 02912, USA;
e-mail: \tt{diana@math.brown.edu}
}
\ Dmitry Fuchs,\thanks{
Department of Mathematics,
University of California, Davis, CA 95616, USA;
e-mail: \tt{fuchs@math.ucdavis.edu}
}
\ and Serge Tabachnikov\thanks{
Department of Mathematics,
Pennsylvania State University, University Park, PA 16802, USA;
e-mail: \tt{tabachni@math.psu.edu}
}
\\
}

\maketitle

\vskip 10 mm

{\hfill To the memory of V. I. Arnold}

\bigskip

\section{Introduction} \label{intro}

The study of billiards in rational polygons and of directional flows on flat surfaces is a fast-growing and fascinating  area of research. A classical construction  reduces the billiard system in a rational polygon -- a polygon whose angles are $\pi$-rational -- to a constant flow  on a flat surface with conical singularities, determined by the billiard polygon. In the most elementary case, the billiard table is a square and the surface is a flat torus obtained from four copies of the square by identifying pairs of parallel sides. We refer  to \cite{HS,MT,Sm,Tab,Vo,Zo} for surveys of rational polygonal billiards and flat surfaces.  

It is well known that the dynamics of a constant flow on a flat torus depends on the direction: if the slope is rational then all the orbits are closed; and if the slope is irrational then all the orbits are uniformly distributed. The same dichotomy holds for the billiard flow in a square. This property is easy to deduce from the fact that a square tiles the plane by reflections in its sides. In the seminal papers \cite{Ve1,Ve2}, W. Veech discovered a large class of polygons and flat surfaces that satisfy the same dynamical dichotomy as the square (it is now called the Veech dichotomy or the Veech alternative). This class includes the isosceles triangle with the angles $(\pi/5,\pi/5,3\pi/5)$ and the regular pentagon; neither tiles the plane by reflection. 

\begin{figure}[hbtp]
\centering
\includegraphics[width=5in]{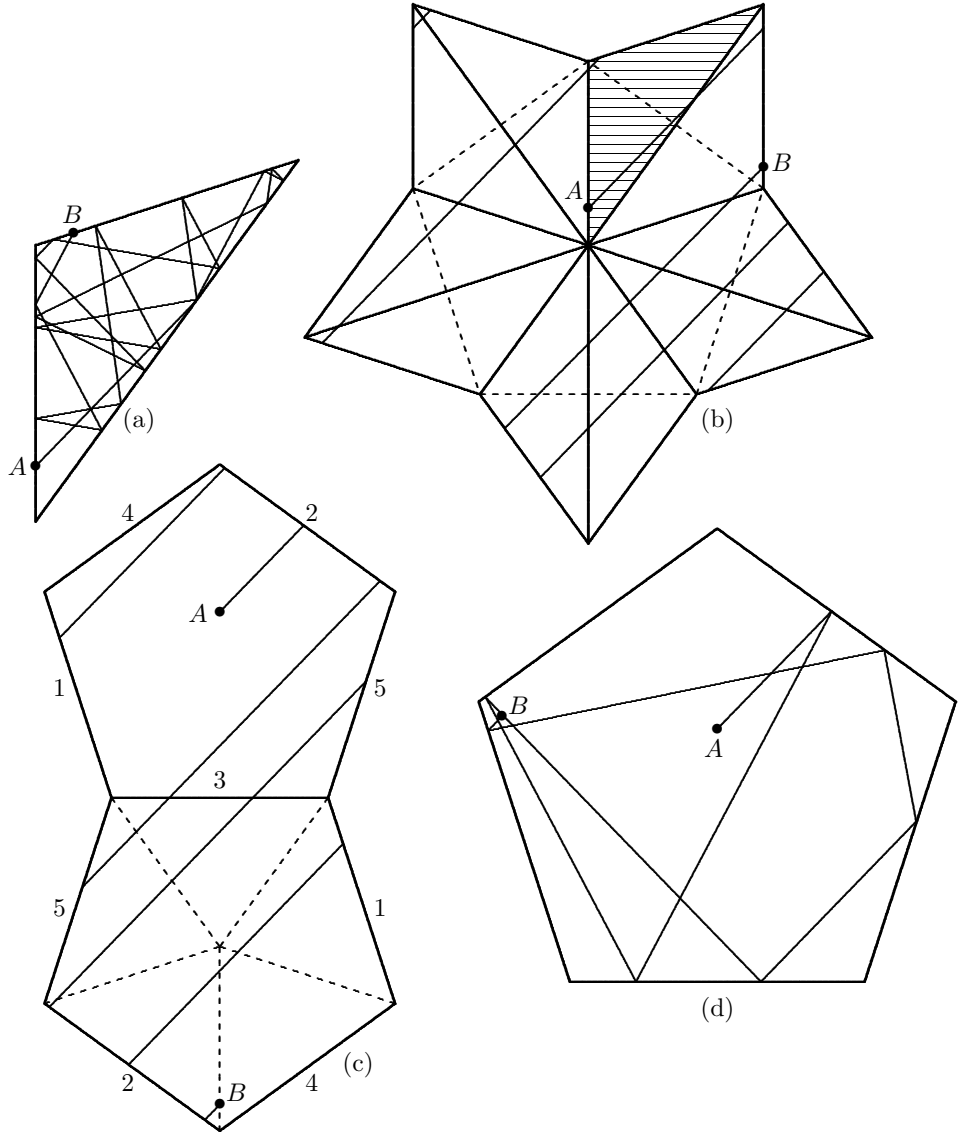} \quad 
\caption{(a) $(\pi/5,\pi/5,3\pi/5)$-triangle with a billiard trajectory; (b) a star-like decagon with the same trajectory; (c) a double pentagon with the same trajectory; (d) a regular pentagon with the same billiard trajectory.}
\label{transition}
\end{figure}

This paper is devoted to a case study of these two polygonal billiards, namely, to a detailed description of periodic billiard trajectories. The flat surface corresponding to the $(\pi/5,\pi/5,3\pi/5)$-triangle is the {\it double pentagon}, an oriented surface of genus 2 constructed from two centrally symmetric copies of the regular pentagon by identifying pairs of parallel sides by parallel translation. The transition between billiards and flat surfaces is shown in Figure \ref{transition}. First, we take 10 copies of the triangle and tile with them a star-like decagon; a billiard trajectory in the triangle (Figure \ref{transition} (a)) becomes a sequence of parallel intervals in the decagon (Figure \ref{transition} (b)) which, in turn, becomes a geodesic on a flat surface of genus 2 obtained from the decagon by attaching parallel sides by means of parallel translation. We rearrange this decagon into a double pentagon by cutting off five triangles by dashed lines and then parallel translations of these triangles so that they form a second (bottom) pentagon (Figure \ref{transition} (c)). The last step shown in Figure \ref{transition} (d) is a transition to a billiard trajectory in a regular pentagon.  The latter is also equivalent to the directional flow on a flat surface; this surface is glued together from ten copies of the regular pentagon and has genus 6. The latter surface admits a 5-fold covering of the former one (see \cite{Ve2}), so a linear periodic trajectory in the double pentagon gives rise to a periodic billiard trajectory in the regular pentagon; the period of the latter may be the same or five times that of the one on the double pentagon; see Section \ref{regpenta}. 

A periodic billiard trajectory in a polygonal billiard is always included into a 1-parameter family of periodic trajectories. If the period is even, one has a strip of parallel trajectories having the same length and the same combinatorial period; if the period is odd, nearby parallel trajectories have the length and the period twice as large. When talking about periodic billiards orbits, we always mean even-periodic ones, and always consider a family of parallel orbits.

Our motivation in this study is two-fold. First, it is an intriguing problem to describe closed geodesics  on the surfaces of regular polyhedra in 3-space, in particular, of a regular dodecahedron \cite{F,FF}. The classifications of closed geodesics on a regular dodecahedron and closed billiard trajectories in a regular pentagon are closely related; we do not  consider the dodecahedron problem here, but we hope to return to it in the near future. Secondly, we were inspired by a recent study by J. Smillie and C. Ulcigrai of the linear trajectories on the flat surface obtained from a regular octagon by identifying the opposite sides \cite{SU1, SU2}. We obtain a number of results for the double pentagon that are analogous to the results in \cite{SU1, SU2}, but we also go further: some of our results are new in the case of the octagon as well. 

\begin{figure}[hbtp]
\centering
\includegraphics[width=5in]{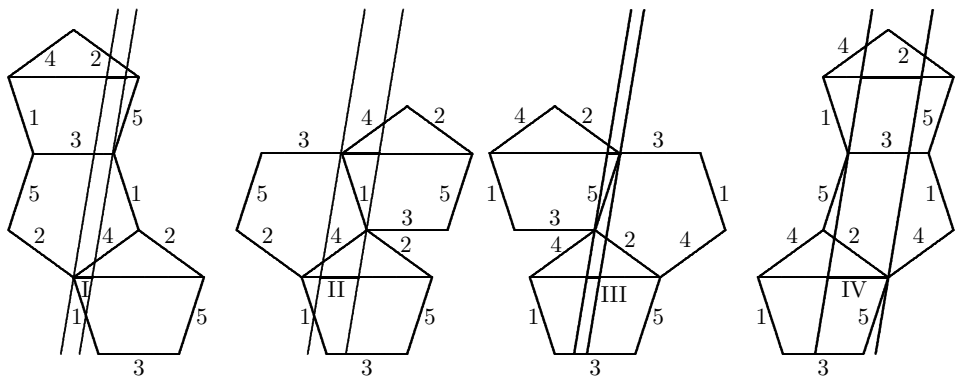}
\caption{Reduction to an interval exchange map}
\label{intex1}
\end{figure}

It is well known that constant flows on flat surfaces are intimately related with interval exchange transformations: one obtains an interval exchange  as a 1-dimensional section of a constant flow and the first return map to this section. In our situation, we have the following equivalent description of the constant flows on the double pentagon depicted in Figure \ref{intex1}. 

\begin{figure}[hbtp]
\centering
\includegraphics[width=1.6in]{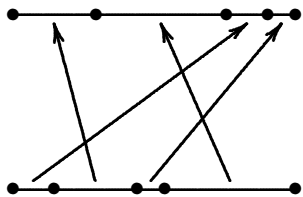}
\caption{Interval exchange map}
\label{intex2}
\end{figure}

Choose a diagonal of the pentagon as a section of a constant flow on the double pentagon and assume that the angle between the trajectory and the diagonal is between $72^\circ$ and $90^\circ$.  An exchange of four intervals arises, permuted as follows:
$\displaystyle{\left(\begin{array} {cccc} {\rm I}&{\rm II}&{\rm III}&{\rm IV}\\ {\rm III}&{\rm I}&{\rm IV}&{\rm II}\end{array}\right)}$, see Figure \ref{intex2}.

If the length of the side of the pentagon is 1, then the length of the diagonal is the golden ratio $\phi=\displaystyle{\frac{1+\sqrt5}2}$, and if we consider the domain of our map as the interval $[0,\phi]$, then the division points are$$\frac12-u(\phi+1),\frac{\phi}2-u,-\frac12-u(\phi-1)+\phi.$$Here $u=\displaystyle{\frac{\sin36^\circ}{\lambda}}$ where $\lambda$ is the slope of the trajectory. Our condition on the angle takes the form $0\le u\le1-\phi/2\approx0.19$. 

We shall use both descriptions, directional flows on the double pentagon and the interval exchange map, interchangeably. 

\bigskip

{\bf Acknowledgments}. D.F. and S.T. are grateful to the Mathematisches Forschungsinstitut Oberwolfach whose hospitality they enjoyed  during their Research in Pairs stay in summer of 2010. D.F. is grateful to IHES for an inspiring atmosphere and excellent working conditions during his visit in summer of 2010. We are grateful to J. Smillie and A. Zorich for interesting discussions, and to R. Schwartz for the graduate topics course on piecewise isometries that he taught in spring of 2010. 

\section{Statements of results} \label{results}

\subsection{Directions}

We start with a remark that the double pentagon has an involution, the central symmetry that exchanges the two copies of the regular pentagon. This involution interchanges linear trajectories having the opposite directions. For this reason, when talking about the direction of a trajectory, we do not make a distinction between the directions $\alpha$ and $\alpha+\pi$. In other words, the set of directions is the real projective line $\RP^1$.

\begin{figure}[hbtp]
\centering
\includegraphics[width=3.2in]{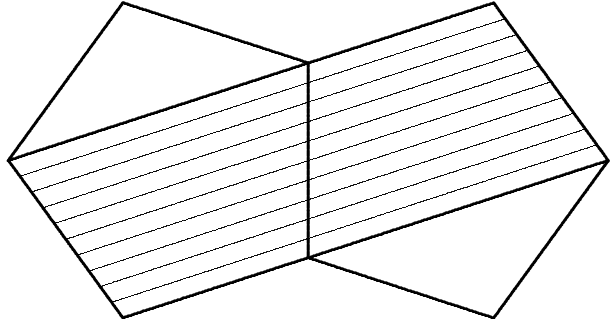}
\caption{Two strips covering the double pentagon: a shadowed wider and longer  strip and a white narrower and shorter strip}
\label{twostrips}
\end{figure}

We identify $\RP^1$ with the circle at infinity (absolute)  of the hyperbolic plane in the Poincar\'e model. The five directions of the diagonals of the pentagon are periodic: the double pentagon decomposes into two strips of parallel periodic trajectories, see Figure \ref{twostrips}. 

\begin{figure}[hbtp]
\centering
\includegraphics[width=4.2in]{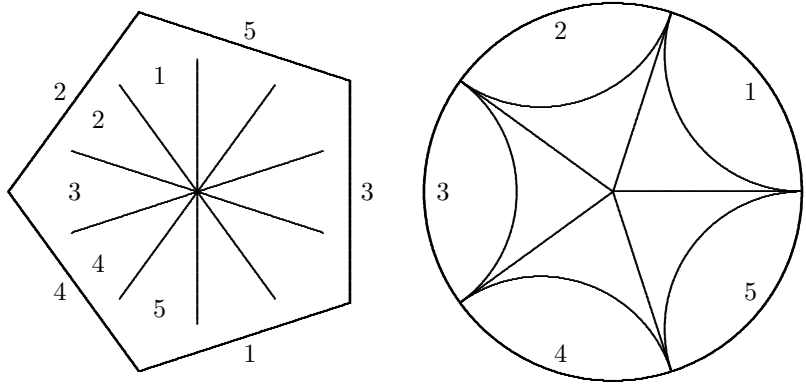}
\caption{Sectors of directions on the absolute}
\label{intervals}
\end{figure}

The directions of the diagonals of the pentagon form the vertices of an ideal regular pentagon in the hyperbolic plane.\footnote{An ideal polygon is regular if the cross-ratio of each consecutive quadruple of its vertices is the same.} We call this ideal pentagon the pentagon of $0$th generation.
The five arcs of $\RP^1$ bounded by the vertices of this ideal pentagon correspond to the  five cones of directions in Figure \ref{intervals}. We consider the 3rd sector as a principal one and focus on periodic directions in this sector.

 We choose an affine coordinate on $\RP^1=\R \cup \infty$ in such a way that the vertices of the ideal regular pentagon of $0$th generation have the coordinates 
$\displaystyle{1-\frac{\phi}2,\frac{\phi}2,\infty,-\frac{\phi}2,\frac{\phi}2-1}$ 
where $\displaystyle{\phi=\frac{1+\sqrt{5}}{2}}$ is the Golden Ratio (this choice is unique, up to a fractional-linear  transformation).

Let $\Gamma$ be the group of isometries of the hyperbolic plane generated by the clockwise rotation $T$ by $72^\circ$ and the reflection $R$ in the vertical side of the pentagon of $0$th generation;  these transformations act on $\RP^1$ by the formulas
$$T(x)=\frac{2\phi x+3-\phi}{2\phi-4x},\ R(x)=\frac1{4\phi^4x}.$$

\begin{figure}[hbtp]
\centering
\includegraphics[width=4.8in]{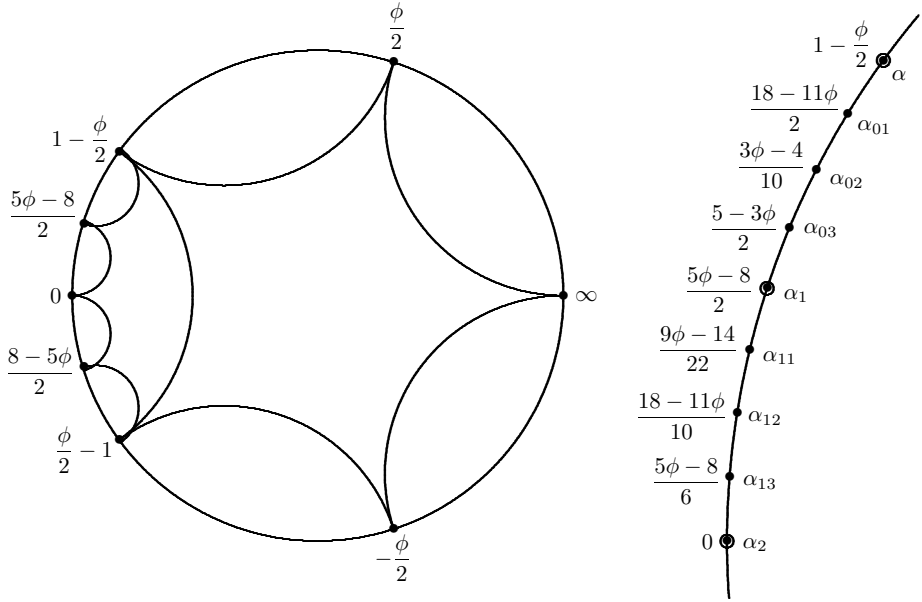}
\caption{Directions of periodic geodesics}
\label{circle}
\end{figure}

The action of $\Gamma$ creates smaller pentagons of the 1st, 2nd, etc., generations, see Figure \ref{circle}. 
We are interested in points to the left of the geodesic  $\left(1-\displaystyle{\phi\over2},\displaystyle{\phi\over2}-1\right)$; on this arc of the absolute the vertices of the pentagon of $k$th generation are denoted by $\alpha$ with $k$ indices as shown in Figure \ref{circle}. For example, the pentagon bounded by the arc $(\alpha,\alpha_1)$ has vertices $\alpha,\alpha_{01},\alpha_{02},\alpha_{03},\alpha_1$, and the pentagon bounded by the arc $(\alpha_{011},\alpha_{012})$ has vertices $\alpha_{011},\alpha_{0111},\alpha_{0112},$ $\alpha_{0113},\alpha_{012}.$ 

The next result (which is not genuinely new) summarizes information about periodic directions. 

\begin{figure}[hbtp]
\centering
\includegraphics[width=3.2in]{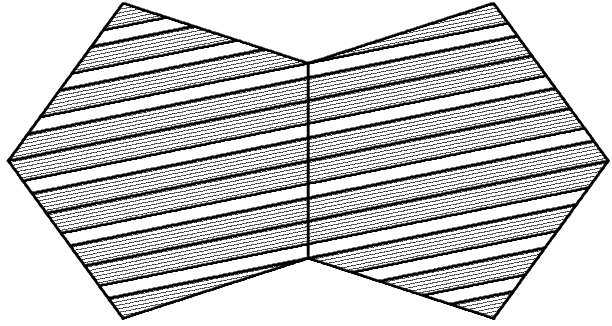}
\caption{Two strips of regular periodic trajectories}
\label{twolongstrips}
\end{figure}

\begin{theorem} \label{dirs}
(i) The set of directions of periodic trajectories (within the 3rd sector) is the set of numbers $\alpha_{n_1n_2\dots n_k}$ with $0\le n_i\le3$ and $n_k\ne0$; \\
(ii) for each of these periodic directions, the double pentagon decomposes into the union of two strips of parallel periodic trajectories (see Figure \ref{twolongstrips});\\
(iii) one has: 
$$\alpha_{n_1n_2\dots n_k}=RT^{m_1}RT^{m_2}\dots RT^{m_k}\alpha$$
where
$$m_i=\left\{\begin{array} {ll} 4-n_i,&\mbox{if}\ i\ \mbox{is even,}\\ n_i+1,&\mbox{if}\ i\ \mbox{is odd and}\ i\ne k,\\ n_i,&\mbox{if}\ i\ \mbox{is odd and}\ i=k;\end{array}\right.$$\\
(iv) the set of periodic directions is $\Q[\phi]\cup\{\infty\}$.
\end{theorem}

\begin{remark}
{\rm
We note that  
$$T^m(x)=-{\phi\over2}-{1\strut\over-\phi-\displaystyle{1\strut\over\displaystyle{\ddots\hskip.6in\atop-\phi-\displaystyle{1\strut\over-\displaystyle{\phi\over2}+x}}}},$$ 
and this allows to transform the formula for $\alpha_{n_1n_2\dots n_k}$ into a continued fraction in the spirit of \cite{Ro1}.}
\end{remark}

\subsection{Periods} \label{periods}

We describe periodic trajectories symbolically in two ways. When dealing with the double pentagon, we label the sides by symbols $1,2,3,4,5$, see Figure \ref{intervals}. Then a periodic linear trajectory has a symbolic orbit, a periodic word in these five symbols, consisting of the labels of the consecutively crossed sides. The (combinatorial) period of the trajectory is the period of this word. When we deal with the interval exchange model, we label the consecutive four intervals by symbols ${\rm I,II,III,IV}$, and then a periodic  trajectory has a symbolic orbit, a periodic word in these four symbols, consisting of the labels of the   visited intervals. 

The two codings translate into each as follows.

\begin{lemma}
The two forms of symbolic orbits correspond to each other according to the rule
$$\{43\}\leftrightarrow\{{\rm I}\};\ \{41\}\leftrightarrow\{{\rm II}\};\ \{25\}\leftrightarrow\{{\rm III}\};\ \{23\}\leftrightarrow\{{\rm IV}\}.$$
\end{lemma}

\proof
(See Figure \ref{intex1}.) We follow upward from the interval I, II, III, IV between the parallel lines and record the labels of the sides we cross.
\proofend

In particular, the period in the Roman numerals is half of the respective period in the Arabic ones.

The next theorem holds for both kinds of periods, ``Roman" and ``Arabic". According to Theorem \ref{dirs}, to every periodic direction, represented by a point on arc of the absolute, there correspond two periods; we denote them by a pair of script and capital letters, such as $a$ and $A$, $a\le A$. 

\begin{theorem} \label{pers}
Let $a,A$ and $b,B$ be the pairs of periods for two points joined by a side of an ideal pentagon of some generation. Then, for the three additional vertices of the pentagon of the next generations, the periods are as shown in the diagram below:

\begin{figure}[hbtp]
\centering
\includegraphics[width=2in]{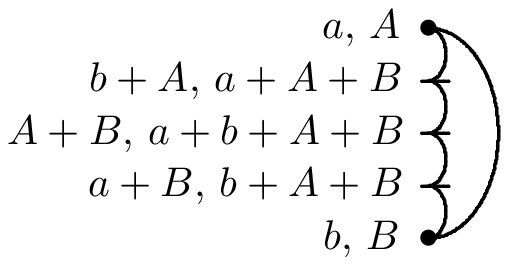}
\label{arcs2}
\end{figure}
\end{theorem}

The (Roman) periods for the points $\pm(1-\phi/2)$ are all equal to 1. Thus, we can compute the periods of all periodic orbits. For example, for the points shown in Figure \ref{circle}, the periods are given in the following table.\bigskip

\centerline{\vbox{\offinterlineskip \hrule \hrule 
\halign{&\vrule#&\strut\hskip4pt#\hskip4pt\cr 
height2pt&\omit&&\omit&&\omit&&\omit&\cr 
&Direction&&Periods&&Direction&&Periods&\cr
height2pt&\omit&&\omit&&\omit&&\omit&\cr 
\noalign{\hrule}
height2pt&\omit&&\omit&&\omit&&\omit&\cr 
&$\quad\alpha$&&\quad1,1&&$\quad\alpha_{11}$&&\quad5,9&\cr
height2pt&\omit&&\omit&&\omit&&\omit&\cr 
\noalign{\hrule}
height2pt&\omit&&\omit&&\omit&&\omit&\cr 
&$\quad\alpha_{01}$&&\quad3,5&&$\quad\alpha_{12}$&&\quad7,11&\cr
height2pt&\omit&&\omit&&\omit&&\omit&\cr 
\noalign{\hrule}
height2pt&\omit&&\omit&&\omit&&\omit&\cr 
&$\quad\alpha_{02}$&&\quad4,7&&$\quad\alpha_{13}$&&\quad6,9&\cr
height2pt&\omit&&\omit&&\omit&&\omit&\cr 
\noalign{\hrule}
height2pt&\omit&&\omit&&\omit&&\omit&\cr 
&$\quad\alpha_{03}$&&\quad4,6&&$\quad\alpha_2$&&\quad2,4&\cr
height2pt&\omit&&\omit&&\omit&&\omit&\cr 
\noalign{\hrule}
height2pt&\omit&&\omit&&\omit&&\omit&\cr 
&$\quad\alpha_1$&&\quad2,3&&&&&\cr
height2pt&\omit&&\omit&&\omit&&\omit&\cr 
\noalign{\hrule}
}}} \bigskip

When proceeding to further generations, the periods grow rapidly. For example, for the direction $\alpha_{123123123}$ the two periods are $3932, 6334$.

There exists an equivalent statement of Theorem \ref{pers}.

\begin{figure}[hbtp]
\centering
\includegraphics[width=3in]{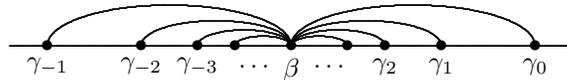}
\caption{Periodic directions connected with a given periodic direction}
\label{arcs}
\end{figure}

\begin{theorem} \label{persprime}
Let $\beta$ be some of the directions $\alpha_{n_1\dots n_k}$, and let $\dots,\gamma_{-2},\gamma_{-1},$ $\gamma_0,\gamma_1,\gamma_2,\dots$ be all points connected by arcs (sides of ideal pentagons) with $\beta$ and ordered as shown in Figure \ref{arcs}. Let $(b,B)$ be the periods corresponding to the direction $\beta$, and let $(c_i,C_i)$ be the periods corresponding to $\gamma_i$ multiplied by $-1$, if $i<0$. Then 
$$\dots,(c_{-2},C_{-2}),(c_{-1},C_{-1}),(c_0,C_0),(c_1,C_1),(c_2,C_2),\dots$$
is an arithmetic sequence with the difference $(B,b+B)$.
\end{theorem}

\proof
We need to check that for any three vertices in a row of an ideal pentagon of any generation, the pairs of periods are as shown in the diagrams below (with some $c,d,C,D$)

\begin{figure}[hbtp]
\centering
\includegraphics[width=4in]{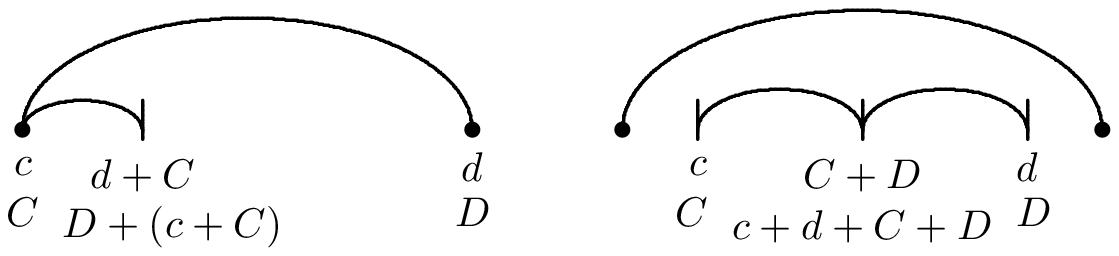}
\label{arcs1}
\end{figure}

\noindent (the left diagram may be reflected in a vertical line; in the right diagram, one of the ends of the upper arc may coincide with one of the three points with periods marked). To do this, it is sufficient to compare these diagrams with corresponding fragments of the diagram in Theorem \ref{pers}.
\proofend

One may present the information on the periods in a different format. Encode the pair of (Roman) periods $(a,A)$ by a single number $a+\phi A \in \Z[\phi]$. For example,  such a period for the points $\pm(1-\phi/2)$ is $1+\phi = \phi^2$. In this notation,  Theorem \ref{pers} asserts that if $U,V\in \Z[\phi]$ are the periods at two consecutive points of $k$th generation (from top to bottom, in Figure \ref{circle}) then the respective periods in the $(k+1)$st generation are
$$
U,\ V+ \phi U,\ \phi U + \phi V,\ U + \phi V,\ V.
$$
This can be restated as follows: the consecutive pair $(U,V)$ is replaced by four consecutive pairs
$$
(U,\ V+ \phi U),\ (V+ \phi U,\ \phi U + \phi V),\ (\phi U + \phi V,\ U + \phi V),\ (U + \phi V,\ V)
$$
obtained from the column vector $(U,V)^{T}$  by the action of following matrices:
$$
X_0=\left(\begin{array}{cc}
1&0\\
\phi&1
\end{array}\right),\ 
X_1=\left(\begin{array}{cc}
\phi&1\\
\phi& \phi
\end{array}\right),\ 
X_2=\left(\begin{array}{cc}
\phi& \phi\\
1& \phi
\end{array}\right),\ 
X_3=\left(\begin{array}{cc}
1& \phi\\
0&1
\end{array}\right).
$$
One obtains the following consequence of Theorem \ref{pers}. 

\begin{theorem} \label{count}
Given a periodic direction $\alpha_{n_1n_2\dots n_k}$, the respective period in $\Z[\phi]$ is the first component of the vector
$$
 X_{n_k} X_{n_{k-1}}\dots X_{n_1} \left(\begin{array}{c}
\phi^2\\
\phi^2
\end{array}\right).
$$
\end{theorem}

\subsection{Symbolic orbits} \label{symborb}

\begin{figure}[hbtp]
\centering
\includegraphics[width=1.6in]{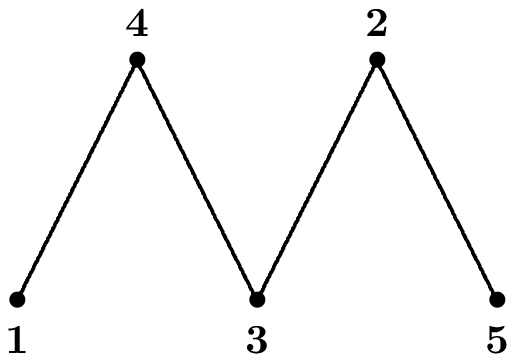}
\caption{Generating graph of symbolic orbits}
\label{graph}
\end{figure}

For every $\alpha_{n_1n_2\dots n_k}\ (k\ge0,\, 0\le n_i\le3,\, n_k\ne0,\ k$ is the number of generation) there are two periodic symbolic orbits, defined up to a cyclic permutation: a ``short" one and a ``long" one (although for $k=0$ they have the same length). Below, we describe an algorithm which creates orbits for the $(k+1)$st generation from the orbits of $k$th generation. In this construction we use symbolic orbits arising from the double pentagon, i.e., the symbols used are $1,2,3,4,5$. The construction is based on the graph depicted in Figure \ref{graph}.

Take a short or long symbolic orbit corresponding to $\alpha_{n_1n_2\dots n_k}$. For example, the short orbit corresponding to $\alpha_{03}$ is $4\ 3\ 2\ 3\ 4\ 1\ 4\ 1$. Shift this (cyclically) by $j=1,2,3,4$; in our example we have:$$\begin{array} {lrrrrrrrr} (\mbox{for}\ j=1)&5&4&3&4&5&2&5&2\\ (\mbox{for}\ j=2)&1&5&4&5&1&3&1&3\\ (\mbox{for}\ j=3)&2&1&5&1&2&4&2&4\\ (\mbox{for}\ j=4)&3&2&1&2&3&5&3&5 \end{array}$$Then, for every of the four sequences, we travel along our graph from each symbol in the sequence to the next one, and insert all the symbols we pass. Our sequences become longer; in our example, we get the following four sequences (we use bold face font for the symbols which come from the old sequence):
\begin{equation} \label{longer}
\begin{array} {cccccccccccccccccccccc} {\bf5}&2&3&{\bf4}&{\bf3}&{\bf4}&3&2&{\bf5}&{\bf2}&{\bf5}&{\bf2}\cr {\bf1}&4&3&2&{\bf5}&2&3&{\bf4}&3&2&{\bf5}&2&3&4&{\bf1}&4&{\bf3}&4&{\bf1}&4&{\bf3}&4\cr {\bf2}&3&4&{\bf1}&4&3&2&{\bf5}&2&3&4&{\bf1}&4&3&{\bf2}&3&{\bf4}&3&{\bf2}&3&{\bf4}&3\cr {\bf3}&{\bf2}&3&4&{\bf1}&4&3&{\bf2}&{\bf3}&2&{\bf5}&2&{\bf3}&2&{\bf5}&2\cr
 \end{array}
\end{equation}

\begin{theorem} \label{symb1}
The four sequences obtained are symbolic orbits at the points$$\alpha_{j-1,3-n_1,3-n_2,\dots,3-n_{k-1},4-n_k},\ j=1,2,3,4,$$short, if the initial orbit is short, and long, if the initial orbit is long.
\end{theorem}

In the example above, the four sequences obtained are short symbolic orbits at the points $\alpha_{031},\alpha_{131},\alpha_{231},\alpha_{331}$.

A symbol of a symbolic orbit is called {\it sandwiched} if the symbols preceding it and following it are the same. Note that all boldface symbols in (\ref{longer}) are sandwiched. The {\it reduction} of a symbolic orbit is obtained by keeping only the sandwiched symbols and deleting all other ones. For example, the reduction of the periodic word 
$5\ 2\ 3\ 4\ 3\ 2\ 3\ 4\ 3\ 2$ is the periodic word $4\ 2\ 4\ 5$. 

The next result is converse to Theorem \ref{symb1}. A similar sandwiching property for octagons was discovered in \cite{SU1}.

\begin{theorem} \label{sand}
Consider a symbolic orbit, short or long, corresponding to a periodic direction  $\alpha_{n_1n_2\dots n_k}$. Reduce this symbolic orbit and shift the reduced word (cyclically) by $4-n_1$. The resulting cyclic word is the symbolic orbit corresponding to the periodic direction  $\alpha_{3-n_2,\dots,3-n_{k-1},4- n_k}$, short or long, respectively. 
\end{theorem}

To translate the symbolic orbits in the double pentagon into the language of the interval exchange map, we need to make the ``inverse change" $\{43\}\to {\rm I},\, \{41\}\to {\rm II},\, \{25\}\to {\rm III},\, \{23\}\to {\rm IV}$. 

\begin{lemma}
Every cyclic symbolic orbit in symbols $\{1,2,3,4,5\}$ can be written as a cyclic word in symbols $\{{\rm I,II,III,IV}\}$.
\end{lemma}

\proof
The pairs in question are precisely the edges of the graph in Figure \ref{graph} oriented downward. A path in this graph has upward and downward edges alternating, hence the sequence of vertices passed can be split into pairs corresponding to downward edges.
\proofend

 Assign to a ``Roman" cyclic symbolic orbit the 4-vector whose components are equal to the number of symbols ${\rm I,II,III,IV}$ in the orbit. We denote the vectors corresponding to the short and the long orbits in the same direction by $(c,d,e,f)$ and $(C,D,E,F)$.

\begin{theorem} \label{relvect}
One has: 
$C=c+e,\, D=f,\, E=c,\, F=d+f.$
\end{theorem}

\subsection{Further experimental results}

We state below two more propositions concerning symbolic orbits. They are confirmed by a huge experimental material, and we hope to give their proofs in forthcoming publications. They can be regarded as symbolic counterparts of Theorems \ref{pers} and \ref{persprime}.

\begin{conjecture} \label{concat}
Let $a,A$ and $b,B$ be two pairs of cyclic symbolic orbits corresponding to two periodic directions  joined by a side of an ideal pentagon of some generation. Then one can cut the cyclic words $a,A,b,B$ into linear ones, concatenate them, and close the words up to cyclic words, so that 
the cyclic symbolic orbits for  the three additional vertices of the pentagon of the next generations (listed in the direction from the first point to the second)  are:
$$\begin{array} {ll} \ bA,&BaA;\\ A B, &bBaA;\\ aB,&AbB.\end{array}$$
\end{conjecture}

\begin{conjecture} \label{symbprime}
Let $\beta$ and $\dots,\gamma_{-2},\gamma_{-1},\gamma_0,\gamma_1,\gamma_2,\dots$ denote the same as in Theorem \ref{persprime}. Then there exist a splitting of the short symbolic orbit corresponding to $\beta$ into two parts, $(c,d)$, and two splittings of the long symbolic orbit corresponding to $\beta$ into two parts: $(a,b)=(a',b')$ such that $a$ and $b'$ have the same beginning and such that the short and long symbolic orbits corresponding to $\gamma_i$ look as shown in the diagram below
\begin{figure}[hbtp]
\centering
\includegraphics[width=5in]{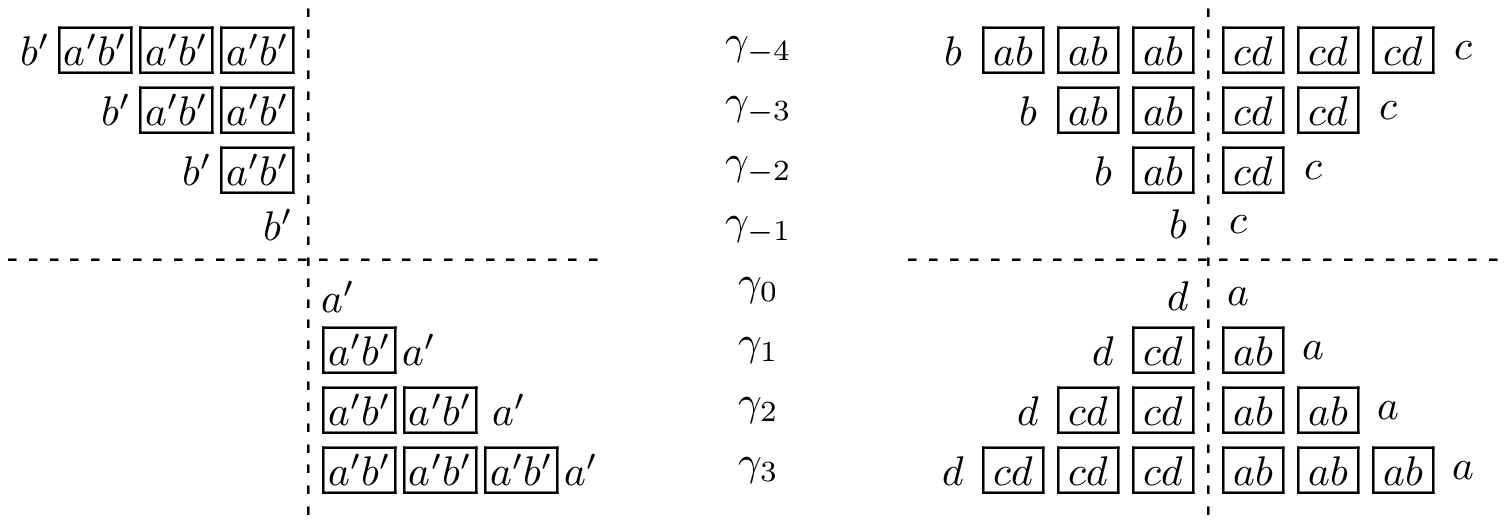}
\label{orbits}
\end{figure}
\end{conjecture}

For example, the symbolic orbits for $\beta=\alpha_{11}$ (in the exchange interval version) are ${\rm IV, I, IV, II, I}$ and ${\rm III, IV, III, IV, II, I,IV, II, I}$. Their splittings from Theorem \ref{symbprime} are

\begin{figure}[hbtp]
\centering
\includegraphics[width=2.8in]{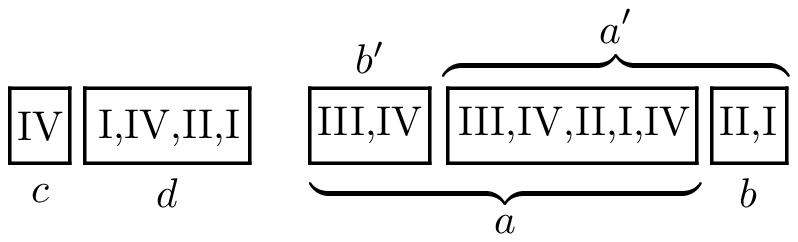}
\label{orbitsfine}
\end{figure}

Then $\gamma_1=\alpha_{111}$ and the symbolic orbits are
$$\underbrace{\rm III, IV, II, I, IV, II, I,}_{a'}\underbrace{\rm III, IV,}_{b'}\underbrace{\rm III, IV, II, I, IV, II, I}_{a'}
$$
and
$$ \underbrace{\rm I,IV,II,I,}_d\underbrace{\rm IV,}_c\underbrace{\rm I,IV,II,I,}_d\underbrace{\rm III,IV,III,IV,II,I,IV,}_a\underbrace{\rm II,I,}_b\underbrace{\rm III,IV,III,IV,II,I,IV}_a$$

\subsection{Lengths and displacement vectors}

A linear trajectory on a flat surface unfolds to a straight line in the plane; a periodic linear orbit  develops to a vector that we call the {\it displacement vector} of the periodic orbit. The displacement vector contains information about the direction and the length of the periodic orbit.

Let the vectors $(c,d,e,f)$ and $(C,D,E,F)$ have the same meaning as before. 

\begin{figure}[hbtp]
\centering
\includegraphics[width=1.9in]{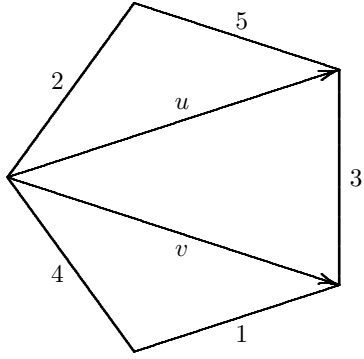}
\caption{Vectors $u$ and $v$}
\label{vectors}
\end{figure}

\begin{theorem} \label{displ}
The displacement vectors of the short and the long periodic orbits are 
$$
(c\phi+e)u+(f\phi+d)v\quad {\rm and}\quad  (C\phi+E)u+(F\phi+D)v
$$
where $u$ and $v$ are the two diagonals bounding the 3rd sector, see Figure \ref{vectors}. The length of the respective short periodic orbit equals
$$m[(c+f)\phi+(d+e)],\ m=\frac{\phi^2\sin36^\circ}{\cos\alpha}$$
where $\alpha$ is the angle between the trajectory and the bisector of the vectors $u$ and $v$; the respective long periodic orbit is $\phi$ times as long.
\end{theorem}

\subsection{Periodic billiard orbits in the regular pentagon} \label{regpenta}

We can relate closed billiard trajectories in the regular pentagon with linear periodic trajectories in the double pentagon.

\begin{theorem} \label{relate}
Let $L$ be the length of closed trajectory of the double pentagon and $(c,d,e,f)$ be the respective 4-component vector.
Then the billiard trajectory in the regular pentagon, starting from the same point in the same direction is also periodic and its length 
is 
$$\left\{\begin{array} {rl} L,&{\rm if}\ (c-f)+2(e-d)\equiv0\bmod5\\ 5L,&{\rm if}\ (c-f)+2(e-d)\not\equiv0\bmod5 \end{array}\right.$$
\end{theorem}

\proof
Consider the left diagram in Figure \ref{intex1}. If we denote the vertices of the bottom pentagon as $A,B,C,D,E$ and then make reflections in the side $BC$ of this pentagon and in the side $AE$ of the reflected pentagon, we will see that the pentagon is rotated by the angle $2\pi/5$ clockwise. 

\begin{figure}[hbtp]
\centering
\includegraphics[width=1.6in]{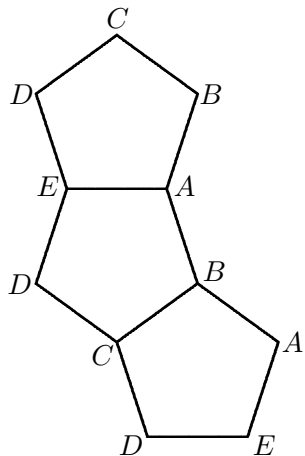}
\caption{Rotation of the pentagon}
\label{rotation}
\end{figure}

Similarly, the compositions of two reflections corresponding to the other three diagrams of Figure \ref{intex1} lead to rotations of the pentagon, respectively, by $4\pi/5$ counterclockwise, by $4\pi/5$ clockwise, and by $2\pi/5$ counterclockwise. If we apply these four transformations, respectively, $c,d,e$, and $f$ times, then the total (clockwise) rotation will be $2((c-f)+2(e-d))\pi/5$. If the latter is not a multiple of $2\pi$, we need to traverse the double pentagon (or interval exchange) periodic orbit 5 times to obtain a closed billiard trajectory in the regular pentagon.
\proofend

\begin{corollary} \label{ratio}
The ratio of the lengths of the long and short closed billiard trajectories of the same directions is always equal to $\phi$.
\end{corollary}

\proof One has:
$(C-F)+2(E-D)=(c+e-d-f)+2(c-f)=3(c-f)+(e-d)\equiv3[(c-f)+2(e-d)]\bmod5.$ \proofend

\begin{remark} \label{dod}
{\rm  A similar fact does not hold for geodesics on a regular dodecahedron: examples show that the ratio between the lengths of parallel geodesics may be different for different directions.}
\end{remark}

\section{Proofs} \label{proofs}

\subsection{Double pentagon as a translation surfaces and its Veech group} \label{transl}

In this section we briefly review basic facts about translation surfaces (our double pentagon is one) and their symmetries. See \cite{HS,MT,Sm,Zo} for a comprehensive exposition. We specify the general theory to the case of the double pentagon and establish Theorem \ref{dirs}. 

A translation surface is a closed surface with conical singularities equipped with an atlas for which the transition functions are parallel translations. 
For our purposes, we can define a translation surface as a surface that is obtained from a finite collection of plane polygons by identifying pairs of parallel sides by translations. 

The group $SL(2,\R)$ acts on plane polygons preserving parallel sides, hence  $SL(2,\R)$ acts on the space of translation surfaces. The Veech group $V_+(X)$ of a translation surface $X$ consists of those $g\in SL(2,\R)$ for which $g(X)$ is equivalent to $X$ as a translation surface. The affine group $Aff_+(X)$ consists of orientation-preserving homeomorphisms of $X$ that are affine in local coordinates; the derivative of such an affine diffeomorphism is constant. The group $V_+(X)$ consists of the derivatives of the transformations in  $Aff_+(X)$.

Following Smillie and Ulcigrai \cite{SU1}, we extend the above described groups to include orientation reversing transformations: $Aff(X)$ consists of all affine diffeomorphisms of $X$, and $V(X)$ of their derivatives. Of course, elements of $Aff(X)$ take periodic trajectories to periodic ones.

A translation surface $X$ is called a Veech surface if $V_+(X)$ is a lattice in $SL(2,\R)$ (a discrete finite co-volume subgroup). For Veech surfaces, the dynamical dichotomy described in the Section \ref{intro} holds. 

In \cite{Ve1},  Veech described the group $V(X)$ for the double pentagon $X$ and proved that it is a Veech surface. This description is crucial for our purposes. 

First of all, symmetries of the regular pentagon provide elements of the group $Aff(X)$. Consider the decomposition of the double pentagon into two horizontal strips in Figure \ref{twostrips}. Choose the direction of the strips as horizontal, and consider the horizontal Dehn twist, a shear map of both strips (a shear map is given by the formula $(x,y)\mapsto (x+cy,y)$) that leaves the horizontal boundaries intact and wraps vertical segments around the strips once. If the side of the pentagon is unit then  $c=2\cot(\pi/5)$. Borrowing  from \cite{SU1}, post-compose the horizontal Dehn twist with reflection in the vertical line; we obtain an affine automorphism $\Phi\in Aff(X)$. The derivative of $\Phi$ is
$$
R=\left(\begin{array}{cc}
-1&-2\cot\displaystyle{\left(\frac{\pi}{5}\right)}\\
0&1
\end{array}\right)
$$
Note that $R$ is an involution, and that $R$ fixes the boundary directions of sector 3 in Figure \ref{intervals}. Note also that $R$ takes the union of sectors $1,2,4,5$ to sector $3$.

The group $V(X)$ is generated by $R$ and the group of symmetries of the regular pentagon $D_5$. We think of elements of the group $V(X)$ as isometries of the hyperbolic plane. Then $R$ is a reflection in the vertical side of the regular ideal pentagon in Figure \ref{intervals}, and $V(X)$ is generated by reflections in the sides of the hyperbolic triangle with angles $(\pi/2,\pi/5,0)$. 

\begin{figure}[hbtp]
\centering
\includegraphics[width=2in]{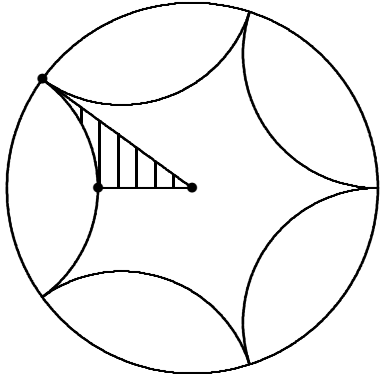}
\caption{$(\pi/2,\pi/5,0)$ triangle}
\label{triangle}
\end{figure}

The directions of periodic linear trajectories on $X$ correspond to the cusps of this group. This set of periodic directions is the set of vertices of the tiling of the hyperbolic plane by ideal pentagons, obtained from the  one in Figure \ref{intervals} by reflections in the sides (the vertices of this ideal pentagon of generation $0$ represent the periodic directions parallel to the sides of the double pentagon).

Thus all periodic directions in sector 3 are obtained by an iterative procedure: apply the rotations $T^m,\ m=1,2,3,4$, to a periodic direction of generation $k$ to obtain a periodic direction in one of the four other sectors, and then apply the reflection $R$ to take this point to a periodic direction in sector 3 of the next generation $k+1$. This is the mechanism described in Theorem \ref{dirs} (i). Since $T$ and $R$ correspond to affine autmorphisms of the double pentagon, they preserve the decomposition into two periodic strips in the periodic direction (statement (ii)). The formulas in statement (iii) are verified directly by induction on $k$. 

Concerning the last statement of Theorem  \ref{dirs}, it is clear that  the set of vertices of the tiling by regular ideal pentagons is a subset of $\Q[\phi]\cup\{\infty\}$. The non-trivial fact that this set coincides with $\Q[\phi]\cup\{\infty\}$ follows from the work of A. Leutbecher \cite{Le}, where this fact is proved for the set of cusp points of the Hecke group $G(2\cos(\pi/5))$, the subgroup of $SL(2,\R)$ generated by the transformations
$z\mapsto -1/z,\ z\mapsto z+ 2\cos(\pi/5)$. One can also deduce this statement from a description of periodic directions for genus 2 translation surfaces with one conical singularity recently obtained  by K. Calta \cite{Ca} and by C. McMullen \cite{Mc}. 

\subsection{Generating symbolic orbits} \label{gensymb}

Following \cite{SU1}, we consider the pairs of consecutive symbols that may appear in the symbolic orbits in the double pentagon. We present the result in the form of a graph whose vertices correspond to the labels of the sides; an oriented edge from vertex $i$ to vertex $j$ means that the pair $ij$ appears in some symbolic orbit. The answers depend on the sector under consideration. The next lemma is straightforward.

\begin{figure}[hbtp]
\centering
\includegraphics[width=4in]{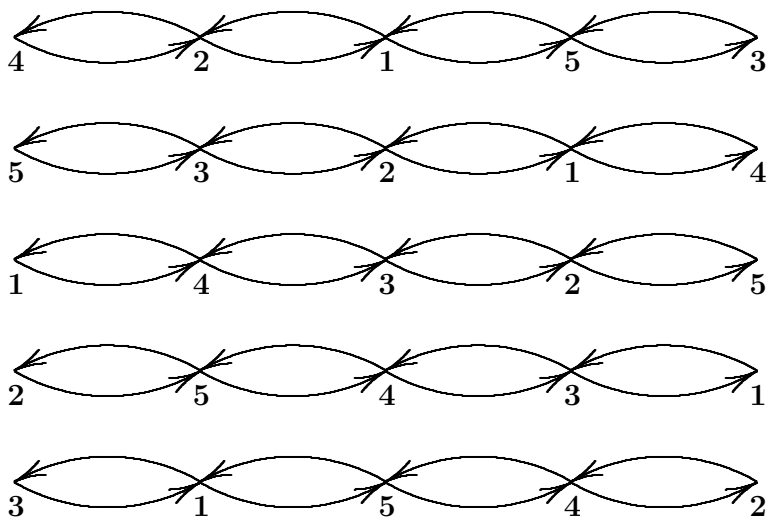}
\caption{Pairs of consecutive symbols in symbolic orbits}
\label{graphs}
\end{figure}

\begin{lemma}  \label{graphlemma}
The graphs describing pair of consecutive symbols in sectors $1,2,3,4,5$, respectively, are depicted in Figure \ref{graphs}.
\end{lemma}

Thus a periodic symbolic trajectory in $i$th sector is a periodic path on $i$th graph in Figure \ref{graphs}. 

\begin{figure}[hbtp]
\centering
\includegraphics[width=4in]{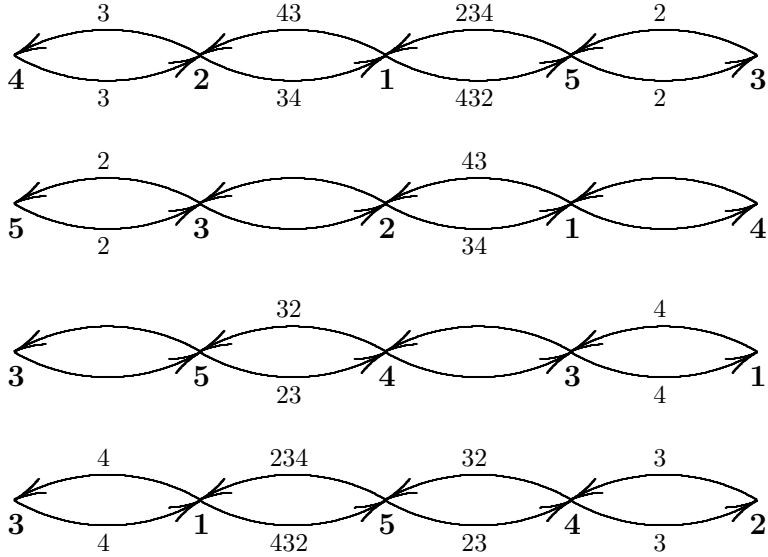}
\caption{Enhanced graphs}
\label{enhgraphs}
\end{figure}

Next we examine the effect of the affine map $\Phi$ on symbolic orbits in sectors $1,2,4,5$; the results are symbolic orbits in sector $3$. We present the result via enhanced graphs with additional symbols written on the edges, see Figure \ref{enhgraphs}. For a symbolic orbit  in sector $i$, presented as a periodic path on $i$th graph, one  traverses this path, inserting the words written on the edges each time the respective edge is passed. We call this the {\it enhancement} of a symbolic  orbit.

\begin{lemma}
Let $w$ be a  symbolic orbit of a periodic trajectory in sector $i$. Apply the affine automorphism $\Phi$ to this trajectory and let $w'$ be its symbolic orbit (which is a periodic path on $3$rd graph in Figure \ref{graphs}). Then $w'$ is the enhancement of $w$.
\end{lemma}

\proof
We consider the case of 1st sector; the other ones are similar.

\begin{figure}[hbtp]
\centering
\includegraphics[width=2.6in]{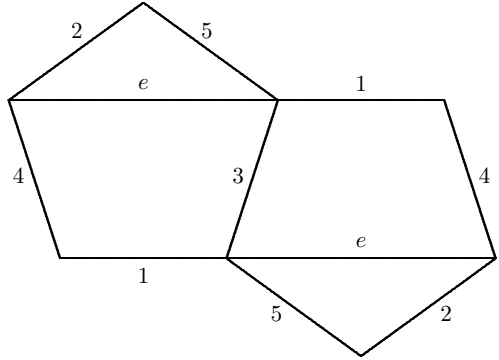}
\caption{Added diagonal}
\label{twist}
\end{figure}

Draw the horizontal diagonal in both copies of the pentagon that make the double pentagon and label this diagonal $e$, see Figure \ref{twist}. We add the symbol $e$ to our alphabet, so the symbolic trajectories will be periodic words in $\{1,2,3,4,5,e\}$. Every trajectory in sector 1 intersects a diagonal marked $e$ between every two consecutive intersections with the sides. That is, we insert $e$ between every two symbols $1,2,3,4,5$, which amount to writing $e$ on each oriented edge of 1st graph in Figure \ref{graphs}.

The shear map affects each segment of a trajectory between the horizontal sides of the two strips into which the double pentagon is split: such a segment is modified by adding one turn around the strip. Symbolically, this is described by the transformations:
$$
e3e\mapsto e4e,\ e4e\mapsto e434e,\ e2e\mapsto e5e,\ e5e\mapsto e525e,\ 1e\mapsto 134e,\ e1\mapsto e431.
$$
A reflection in the vertical line is, symbolically, the involution
$$
2\leftrightarrow 5,\ 3\leftrightarrow 4,\ 1\leftrightarrow 1,
$$
so the symbolic action of $\Phi$ is as follows:
$$
e3e\mapsto e3e,\ e4e\mapsto e343e,\ e2e\mapsto e2e,\ e5e\mapsto e252e,\ 1e\mapsto 143e,\ e1\mapsto e341.
$$

\begin{figure}[hbtp]
\centering
\includegraphics[width=4in]{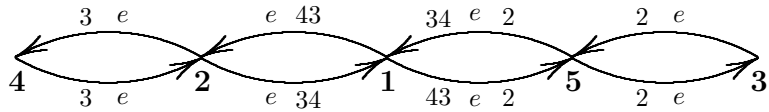}
\caption{Enhanced 1st graph}
\label{graph1}
\end{figure}

This is achieved by enhancement of the graph as shown in Figure \ref{graph1}. It remains to erase the label $e$, and we obtain the 1st enhanced graph in Figure \ref{enhgraphs}.
\proofend

We now observe that all four enhanced graphs in Figure \ref{enhgraphs} are obtained from the graph in Figure \ref{graph} as described in Section \ref{symborb}. 

Let us summarize: given a symbolic periodic orbit in sector 3, apply cyclic shifts to obtain a symbolic trajectories in sectors 1,2,4,5, and enhance these symbolic trajectories by traversing them on the graph in Figure \ref{graph}. This results in symbolic trajectories in sector 3 of the next generation. This is the generation process described in Theorem \ref{symb1}. 

We  also deduce Theorem \ref{sand}. Observe that the sandwiching property holds for all four enhanced graphs in Figure \ref{enhgraphs}: the reduction of an enhanced symbolic orbit is the original orbit. 

For example, consider the path in the 1st graph in Figure \ref{enhgraphs} that visits the vertices $ \dots 3\ 5\ 1\ 2\ 4\dots$ Its enhancement is
$$
\dots 2\ {\bf 3}\ 2\ {\bf 5}\ 2\ 3\ 4\ {\bf 1}\ 4\ 3\ {\bf 2}\ 3\ {\bf 4}\ 3 \dots
$$
and the reduction is again  $ \dots 3\ 5\ 1\ 2\ 4 \dots$.

\subsection{Further proofs} \label{furthpfs}

A symbolic periodic trajectory in sector 3 is a periodic path on graph in Figure \ref{graph}. Assign to such a path the  vector $(c,d,e,f)$ whose components are the number of times per period that the edges $43,\ 41,\ 25,\ 23$ were traversed (in either direction). Thus two such vectors are assigned to every periodic direction; the symbolic period is the sum of the components.
For example, the respective pairs of vectors for the first five periodic directions
$$
1-\frac{\phi}{2},\ \frac{5\phi-8}{2},\ 0,\ \frac{8-5\phi}{2},\ \frac{\phi}{2}-1
$$
in Figure \ref{circle} are:
$$
\begin{array}{ccccccccccccc}
(0&0&1&0),& (1&0&0&0)\\
(1&1&0&0),&(1&0&1&1)\\
(1&0&0&1),&(1&1&1&1)\\
(0&0&1&1),&(1&1&0&1)\\
(0&1&0&0),&(0&0&0&1).
\end{array}
$$ 
Note that if $(a,A)$ and $(b,B)$ are the first and the last pairs of these vectors, then the three pairs in-between are
\begin{equation} \label{rel}
(b+A,a+A+B),\ (A+B,a+b+A+B),\ (a+B,b+A+B).
\end{equation}

Theorem \ref{symb1}, or the enhanced graphs in Figure \ref{enhgraphs}, tell us how the 4-component vector changes under the generation process: the cyclic shift by $i\in\{1,2,3,4\}$ and enhancement. In each case, the result is a linear transformation, depending on $i$. These transformations, denoted by $L_i$, are as follows:
$$
L_1(c,d,e,f)=(c+e+f,e,c+d,c), L_2(c,d,e,f)=(c+d+e+f,c+d,c+f,c+e+f),
$$
$$
L_3(c,d,e,f)=(c+d+f,c+f,e+f,c+d+e+f), L_4(c,d,e,f)=(f,e+f,d,c+d+f).
$$
It follows that the linear relation (\ref{rel}) is inherited by the consecutive quintuples of pairs of vectors of each next generation. Applying the linear function
$$
p(c,d,e,f)=c+d+e+f
$$
to relation (\ref{rel}), we obtain the same relation for periods, that is, the statement of Theorem \ref{pers}.

Theorem \ref{count} is a reformulation of Theorem \ref{pers}, as explained in the paragraph that precedes its formulation in Section \ref{periods}; we do not dwell on its proof.

We can deduce Theorem \ref{relvect} from relation (\ref{rel}). The statement of Theorem \ref{relvect} is that the 4-component vectors $(a,A)$ corresponding to a periodic direction satisfy the linear relation $A=M(a)$ where 
$$
M=\left(
\begin{array}{cccc}
1&0&1&0\\
0&0&0&1\\
1&0&0&0\\
0&1&0&1
\end{array}
\right).
$$
Note that $M^2=M+I$.

Assuming that $A=M(a), B=M(b)$, we want to deduce the same relations for the vectors in (\ref{rel}). Indeed,
$$
a+A+B=a+M(a)+M(b)=M^2(a)+M(b)=M(A+b),
$$
$$
a+b+A+B=a+b+M(a)+M(b)=M^2(a)+M^2(b)=M(A+B),
$$
$$
b+A+B=b+M(a)+M(b)=M(a)+M^2(b)=M(a+B),
$$
and Theorem \ref{relvect} follows.

\begin{figure}[hbtp]
\centering
\includegraphics[width=4in]{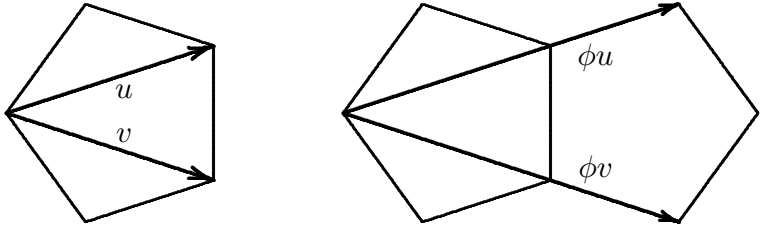}
\caption{Displacement vectors}
\label{displacement}
\end{figure}

Let us also prove Theorem \ref{displ}. The components of the vector $(c,d,e,f)$ are the numbers of symbols ${\rm I,II,III,IV}$ in a periodic symbolic orbit. When the respective periodic trajectory is unfolded in the plane, each symbol ${\rm III}=25$ corresponds to the vector $u$, the symbol ${\rm I}=43$ to the vector $\phi u$, the symbol ${\rm II}=41$ to the vector $v$, and the symbol ${\rm IV}=23$ to the vector $\phi v$, see Figure \ref{displacement}. The displacement vector is then $(c\phi+e)v+(f\phi+d)u$, as claimed. 

For a regular pentagon with the side length 1, the distance from a vertex to the nearest diagonal not passing through this vertex equals $\sin36^\circ$; the distance between a side and the diagonal parallel to this side equals $\sin72^\circ=\phi\sin36^\circ$. Hence, the displacement corresponding to symbols II and III, measured in the directions of the bisector between $u$ and $v$ equals $(1+\phi)\sin 36^\circ=\phi^2\sin36^\circ$, and the displacement, corresponding to the symbols I and IV, measured in the same direction, is $(1+2\phi)\sin36^\circ=\phi^3\sin36^\circ$. The full displacement is $((c+f)\phi^3+(d+e)\phi^2)\sin36^\circ$. To find the length of the actual orbit, we need to divide the result by the cosine of the angle between the direction and the bisector, and we obtain the formula stated in Theorem \ref{displ}.

The two displacement vectors are obtained from the ones, corresponding to the simplest periodic orbits (in the direction $u$), by a linear transformation. The displacement vectors for the latter are $u$ and $\phi u$. It follows that the two vectors are also proportional with coefficient $\phi$:
$$
(C\phi+E)v+(F\phi+D)u = \phi ((c\phi+e)v+(f\phi+d)u).
$$
This again implies the relations of Theorem \ref{relvect}.

\section{Expected further results}

Results similar to the ones of these paper may hold for billiards in all regular polygons and linear flows on double odd-gons or  regular even-gons (in both cases, the parallel sides are identified  by translations). As we mentioned earlier, the case of regular octagon  was studied recently by J. Smillie and C. Ulcigrai (\cite{SU1,SU2}). In addition to the results of these works, we can provide a description of periods of periodic orbits. 

\begin{figure}[hbtp]
\centering
\includegraphics[width=2.4in]{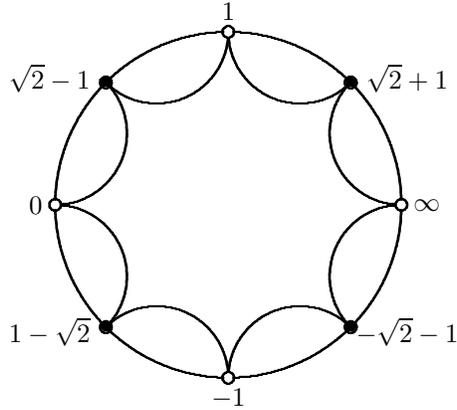}
\caption{Directions of periodic geodesics for an octagon}
\label{octagon}
\end{figure}

A billiard trajectory in a regular octagon is closed if the slope of the trajectory  belongs to $\Q[\sqrt2]$ (we assume that some side is horizontal). These directions can be arranged on the projective line -- the absolute of the hyperbolic plane -- in the following way. First we consider an ideal regular octagon; its vertices, in Figure \ref{octagon}, are marked by the slopes.

As it is seen in Figure \ref{octagon}, the vertices are marked, alternatively, black and white. Then we construct the first generation of octagons by reflection in the sides of the big octagon, then the next generation, and so on. For example, the vertices of the octagon of the first generation between the points $\sqrt2-1$ and 0 are 
$$\sqrt2-1,\, \frac{2\sqrt2-1}7,\, \frac{3-\sqrt2}7,\, \frac{\sqrt2-1}2,\, \frac{3\sqrt2-1}{17},\, 3-2\sqrt2,\frac{\sqrt2-1}3,\, 0,$$ and the vertices of the octagon of the second generation between  $\sqrt2-1$ and $\displaystyle\frac{2\sqrt2-1}7$ are 
$$\sqrt2-1,\frac{3\sqrt2-2}7,\frac{3\sqrt2+1}{17},\frac{10\sqrt2-9}{17},\frac{15\sqrt2-9}{41}, \frac{1-\sqrt2}2,\frac{11\sqrt2-9}{23},\frac{2\sqrt2-1}7.$$
 All vertices of all octagons are marked, alternatively, by white or black dots, and the reflections, as well as rotations by $90^\circ$, preserve these marking.

\begin{figure}[hbtp]
\centering
\includegraphics[width=2.6in]{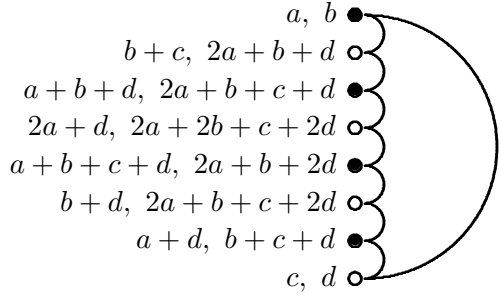}
\caption{Periods in the octagon}
\label{octaper}
\end{figure}

For each of these points, there are two periods. For example, for $\sqrt2-1$, the periods are 1 and 2; for 0, the periods are 2 and 2. A theorem similar to Theorem \ref{pers} states that for the eight vertices of any generation, the pairs of periods look as in Figure \ref{octaper} (one can flip it upside down).\footnote{J. Smillie gave a proof of this result after we informed him about it.}

There also holds a statement similar to Theorem \ref{persprime}, but the difference of the arithmetic sequence depends on the marking of the point $\beta$: if $\beta$ is a white dot, then the difference is $(B,b+B)$, but if is is black, then the difference is $(B,2b+B)$. All other results of Section \ref{results} also have their analogs for octagon.

The results are different for regular heptagons. For those, we have a sequence of generations of ideal regular heptagons, and for every vertex of every  heptagon there are three periods (this 3 is the genus of a flat surface obtained from the double heptagon). If for a heptagon of some generation the triples of periods corresponding to two vertices joined by an arc are $(a_1,a_2,a_3), a_1\le a_2\le a_3$ and $(b_1,b_2,b_3), b_1\le b_2\le b_3$, then the periods at the five intermediate vertices are: 

$$\begin{array} {cl} \bullet&a_1,a_2,a_3\\ -&a_2+b_1, a_1+b_2+a_3, a_2+b_3+a_3\\ -&a_3+b_2, a_2+a_3+b_1+b_3, a_1+a_2+a_3+b_2+b_3\\ -&a_3+b_3, a_2+a_3+b_2+b_3, a_1+a_2+a_3+b_1+b_2+b_3\\ -&a_2+b_3, a_1+a_3+b_2+b_3, a_2+a_3+b_1+b_2+b_3\\ -&a_1+b_2, a_2+b_1+b_3, a_3+b_2+b_3\\ \bullet&b_1, b_2, b_3\end{array}$$

A statement similar to Theorem \ref{persprime} in the case of a regular heptagon looks very simple. If the triple of periods for a point $\beta$ is $(a,b,c)$, then the triples of periods of points $\gamma_i$ (with an appropriate sign change) form an arithmetic sequence with the difference $(b, a+c, b+c)$.

We conjecture that similar results hold, at least, for all regular polygons.

 \end{document}